\documentclass[a4paper]{article}
\usepackage{amsmath,amssymb,amsthm,amscd,euscript}
\DeclareMathOperator{\ord}{ord}

 \DeclareMathOperator{\im}{Im}
 
\DeclareMathOperator{\SO}{SO} \DeclareMathOperator{\Or}{O}
\DeclareMathOperator{\SL}{SL} \DeclareMathOperator{\diag}{diag}
 
\DeclareMathOperator{\Kil}{Kil}

\DeclareMathOperator{\const}{const}

\DeclareMathOperator{\pt}{pt}

\def\pd#1#2{\frac{\partial{#1}}{\partial{#2}}}
\def\pd1#1{\frac{\partial}{\partial#1}}
\def\d1#1{\frac{d}{d#1}}

\def\text#1{\mbox{#1}}

\newtheorem{Def}{Definition}\newtheorem{theore}{Theorem}\newtheorem{Lem}{Lemma}

\allowdisplaybreaks

\numberwithin{theore}{section} \numberwithin{equation}{section}
\numberwithin{proposit}{section}
\numberwithin{Def}{section}\numberwithin{Lem}{section}

\begin{document}

\author{Alexey V.~Shchepetilov\footnote{Department of Physics, Moscow State
University, 119992 Moscow, Russia, e-mail address:
quant@phys.msu.su}}
\title{Nonintegrability of the two-body problem in constant curvature spaces\footnote{\copyright
2006 IOP Publishing Ltd, J. Phys. A: Math. Gen. V.~39 (2006),
pp.~5787-5806}}
\date{}\maketitle
\begin{abstract}
We consider the reduced two-body problem with the Newton and the
oscillator potentials on the sphere ${\bf S}^{2}$ and the
hyperbolic plane ${\bf H}^{2}$. For both types of interaction we
prove the nonexistence of an additional meromorphic integral for
the complexified dynamic systems. \vskip 20pt

\end{abstract}

\section{Introduction}\label{Introduction}\markright{\ref{Introduction}
Introduction}

The study of mechanics on constant curvature spaces begun in the
nineteenth century after the rise of noneuclidean geometry
\cite{Lip3} -- \cite{Neumann}. Similarly to the Euclidean case in
constant curvature simply connected spaces (the sphere
$\mathbf{S}^{n}$ and the hyperbolic space $\mathbf{H}^{n}$) there
are two exceptional central potentials $V_{N}$ and $V_{o}$ (below
Newton and oscillator potentials). They have some nice properties,
which can be grounds for their definitions.

Both these potentials make all bounded trajectories of a one-body
problem closed \cite{Lib03}. Moreover, these trajectories (bounded
and unbounded) are conics \cite{Lib02}, which can be naturally
defined in constant curvature spaces \cite{Story} -- \cite{Klein}.

The one-body motion in the Newton potential satisfies to the
analogues of the three Kepler laws \cite{Kil2}, \cite{Neumann},
\cite{Lib02}. This potential is also the fundamental solution of
the Laplace equation. The corresponding force in the hyperbolic
space was already proposed by N.~Lobachevski (in 1835-38)
\cite{Lob} and J.~Bolyai (between 1848 and 1851) \cite{Bol} as the
value $F(\rho)$ which is inverse to the area of the sphere in
$\mathbf{H}^{3}$ of radius $\rho$ with an attractive body in the
center.

These results can be considered as predecessors of general
relativity. After the rise of this theory the above-mentioned
papers were almost completely forgotten.

Similar models attracted attention later from the point of view of
quantum mechanics and the theory of integrable dynamical systems.
This leads to the rediscovery of results described above in many
papers, sometimes with partial improvements, see for example
\cite{Nishino} and \cite{Higgs}. Note however that almost
forgotten results of W.~Killing and H.~Liebmann were described in
the survey \cite{DomZitt}.

Corresponding quantum mechanical problems in constant curvature
spaces were studied in \cite{Schr} -- \cite{Inf3} and other
papers.

The two-body problem with a central interaction in constant
curvature spaces $\mathbf{S}^{n}$ and $\mathbf{H}^{n}$
considerably differs from its Euclidean analogue. The variable
separation for the latter problem is trivial, while for the former
one no central potentials are known that admit a variable
separation.

The two-body problem with a central interaction in constant
curvature spaces was considered for the first time in
\cite{Shch98}. In Euclidean space this problem is reduced to a
one-body problem in a central potential after separating the
center of mass motion. Due to the absence of Galilei
transformations the situation for the constant curvature spaces is
different. The two-body problem is invariant with respect to the
isometry group, but for non-Euclidean space this group is not wide
enough to imply the integrability of this problem in any sense.

The natural problem of finding central potentials corresponding to
integrable two-body problems is far from its solution now. This
can be explained by the fact that existing methods of the theory
of integrable and nonintegrable dynamical systems do not work in
the presence of a functional parameter.

As a limiting case of a two-body problem in constant curvature
spaces, one can consider the restricted two-body problem: the
"heavy" body moves with a constant velocity along a geodesic,
while the "light" one moves in a potential of a "heavy" body.

The nonintegrability of this problem with the potential $V_{N}$
and $V_{o}$ on the sphere $\mathbf{S}^{2}$ was proved in
\cite{Zi}, \cite{Zi2} in the class of meromorphic functions.
Similar results with smaller restrictions, valid also for the
restricted two-body problem on the hyperbolic plane
$\mathbf{H}^{2}$, were obtained in \cite{MacPrz}.

Here, we prove the nonexistence of an additional meromorphic first
integral for the restricted two-body problem on the spaces
$\mathbf{S}^{2}$ and $\mathbf{H}^{2}$ using the Morales-Ramis
theory \cite{MorRu}.

\section{Reduced two-body problems}\label{RTBP}
\markright{\ref{RTBP} Reduced two-body problems}

Note that the classical two-body problem on $\mathbf{S}^{n}$ and
$\mathbf{H}^{n}$ reaches its full generality at $n=3$
\cite{Shch98}. Its Hamiltonian reduction to the system with two
degrees of freedom was carried out in \cite{Shch98} ($n=3$) and in
\cite{Sh004} ($n=2$) by explicit coordinate calculations. A more
conceptual approach to this reduction was derived in \cite{Sh001}.

Here we shall use the following description of the reduced
dynamical systems for $n=2$, combining approaches from
\cite{Shch98} and \cite{Sh004} -- \cite{Shchep023}.

\subsection{The reduced two-body problem on the sphere $\mathbf{S}^{2}$}

Let ${\bf S}^{2}$ be the sphere of the radius $R$ with the
standard metric. The configuration space for the two-body problem
on ${\bf S}^{2}$ is $Q=\left({\bf S}^{2}\times{\bf
S}^{2}\right)\backslash\diag$. Let $Q_{op}\simeq{\bf S}^{2}$ be a
subset of $Q$, consisting of pairs of opposite points. The phase
space $T^{*}Q$ can be represented as
$$
T^{*}Q=\left(T^{*}I\times
T^{*}\SO(3)\right)\cup\widetilde{T}^{*}Q_{op},
$$
where $I=(0,\pi R)$ and $\widetilde{T}^{*}Q_{op}$ is the
restriction of the cotangent bundle $T^{*}(Q\times Q)$ onto
$Q_{op}$.

The space $\widetilde{T}^{*}Q_{op}$ is the submanifold in $T^{*}Q$
of the codimension $2$; therefore a typical trajectory does not
intersect it. Below we consider only such trajectories.

The group $\SO(3)$ acts by symplectomorphisms on the second factor
of the product
$$
M:=T^{*}I\times T^{*}\SO(3),
$$
endowed with the standard symplectic structure of a cotangent
bundle. Therefore, the reduced phase space for $M$ have the form
\cite{Ar2}
$$
\widetilde{M}=T^{*}I\times\mathcal{O},
$$
where $\mathcal{O}$ is a $\SO(3)$-orbit w.r.t.\ the coadjoint
action in the space $\mathfrak{so}^{*}(3)$ dual to the Lie algebra
$\mathfrak{so}(3)$. The orbit $\mathcal{O}$ is endowed with the
Kirillov symplectic form.

The Killing form on the Lie algebra $\mathfrak{so}(3)$ generates
its natural identification with the dual space
$\mathfrak{so}^{*}(3)$ and makes both these spaces Euclidean. The
coadjoint orbits in $\mathfrak{so}^{*}(3)$ are standart spheres in
the Euclidean space ${\bf E}^{3}$ with the common center
$0\in\mathfrak{so}^{*}(3)$ and the Kirillov symplectic form on
them coincides with area forms, generated by the Euclidean
structure.

The reduced Hamiltonian function on $\widetilde{M}$ for the two
body problem is
\begin{equation}\label{h1}
h_{s}=\frac{(1+r^{2})^{2}}{8mR^{2}}\left(p_{r}^{2}+\frac{p_{2}^{2}}{r^{2}}\right)+
\frac{(1+r^{2})p_{r}p_{0}+\gamma^{2}}{2m_{1}R^{2}}-\frac{p_{2}^{2}}{m_{1}R^{2}}+
\frac{(1-r^{2})p_{1}p_{2}}{2m_{1}R^{2}r}+V(r).
\end{equation}
Here $p_{i}$ are orthogonal coordinates in $\mathfrak{so}^{*}(3)$;
$m:=m_{1}m_{2}/(m_{1}+m_{2})$ is the reduced mass for bodies
masses $m_{1},m_{2}$; $r:=\tan\left(\rho/(2R)\right)$; $\rho$ is
the distance between the bodies, $p_{r}$ is the momentum
corresponding to the coordinate $r$ and the orbit
$\mathcal{O}\equiv\mathcal{O}_{\gamma}$ is defined by the equation
$$
p_{0}^{2}+p_{1}^{2}+p_{2}^{2}=\gamma^{2},\;\gamma\geqslant0.
$$

The Poisson brackets for variables $r,p_{r},p_{0},p_{1},p_{2}$ are
as follows
\begin{align}\begin{split}
\{r,p_{r}\}=1,\,\{p_{0},p_{1}\}=-p_{2},\,\{p_{1},p_{2}\}=-p_{0},\\
\{p_{2},p_{0}\}=-p_{1},\,\{r,p_{i}\}=0,\,\{p_{r},p_{i}\}=0,\,i=0,1,2
\end{split}\label{PoissonS}
\end{align}
and the evolution of any smooth function
$f=f(r,p_{r},p_{0},p_{1},p_{2})$ is defined by the equation
\begin{equation}\label{evolution}
\frac{df}{dt}=\{f,h\}.
\end{equation}

At every fixed moment of time, the momentum $p_{2}$ corresponds to
the rotation of the second body around the first one, the momentum
$p_{0}$ corresponds to the motion of bodies along the geodesic
connecting them and the momentum $p_{1}$ corresponds to the motion
of the system in the direction, normal to this geodesic.

In the exceptional case $\gamma=0$, it holds
$\mathcal{O}_{0}=\pt$; the reduced system has only one degree of
freedom and it corresponds to the motion of the bodies along a
common geodesic with the null value of the total momentum.

Let $\gamma>0$ and $p_{2},\varphi$ be cylinder coordinates on
$\mathcal{O}_{\gamma}$ such that
$$
p_{0}=\sqrt{\gamma^{2}-p_{2}^{2}}\sin\varphi,\,p_{1}=\sqrt{\gamma^{2}-p_{2}^{2}}\cos\varphi.
$$
Then it holds $\varphi=\arctan\left(p_{0}/p_{1}\right)$ and
$\{\varphi,p_{2}\}=1$. Thus $p_{\varphi}:=p_{2},\varphi$ are
canonical coordinates on $\mathcal{O}_{\gamma}$ (with
singularities at the points $p_{2}=\pm\gamma$) and the Hamiltonian
function (\ref{h1}) can be written in the form
\begin{align}\label{h2}
h_{s}&=\frac{(1+r^{2})^{2}}{8mR^{2}}\left(p_{r}^{2}+\frac{p_{\varphi}^{2}}{r^{2}}\right)-
\frac{p_{\varphi}^{2}}{m_{1}R^{2}}
+\frac{\sqrt{\gamma^{2}-p_{\varphi}^{2}}}{2m_{1}R^{2}}
\left(p_{\varphi}\frac{1-r^{2}}r\cos\varphi
+(1+r^{2})p_{r}\sin\varphi\right)\nonumber\\&+V(r)+\frac{\gamma^{2}}{2m_{1}R^{2}},
\end{align}
coinciding with the formula (7) in \cite{Sh004} (cf.\ also the
formula (18) in \cite{Shch98} for $\alpha=\gamma,\beta=0$).

The Hamiltonian function (\ref{h1}) can be represented in another
form after the substitution of the variables $r,p_{r}$ by a new
pair of canonical ones $\theta,p_{\theta}$ such that
$$
\theta=\frac{\rho}{R},\,p_{\theta}=\frac12(1+r^{2})p_{r},\,\{\theta,p_{\theta}\}=1.
$$
This substitution leads to the following expression:
\begin{align}\begin{split}
h_{s}&=\frac{1}{2mR^{2}}\left(p_{\theta}^{2}+\frac{p_{2}^{2}}{\sin^{2}\theta}\right)+
\frac{p_{\theta}p_{0}}{m_{1}R^{2}}-\frac{p_{2}^{2}}{m_{1}R^{2}}+
\frac{p_{1}p_{2}}{m_{1}R^{2}}\cot\theta+\frac{\gamma^{2}}{2m_{1}R^{2}}
+V(\theta)\\&=h_{s,1}+h_{s,2}+\frac{\gamma^{2}}{2m_{1}R^{2}}+V(\theta),
\end{split}\label{hh3}
\end{align}
where
$$
h_{s,1}:=\frac{1}{2mR^{2}}\left(p_{\theta}^{2}+\frac{p_{2}^{2}}{\sin^{2}\theta}\right),\;
h_{s,2}:=\frac{p_{\theta}p_{0}}{m_{1}R^{2}}-\frac{p_{2}^{2}}{m_{1}R^{2}}+
\frac{p_{1}p_{2}}{m_{1}R^{2}}\cot\theta.
$$
Below we shall use expression (\ref{hh3}) for the Hamiltonian
function $h_{s}$ though all claims can be reformulated for other
coordinate systems.

The Hamiltonian function $h_{s}$ corresponds to trivially
integrable Hamiltonian systems in the following three cases.
\begin{enumerate}
\item The free motion of bodies: $V(\theta)=0$. One can easily
verify that $\{h_{s,1},h_{s,2}\}=0$ and an additional integral in
this case is $h_{s,1}$ or $h_{s,2}$.
\item $m_{1}=\infty$. This case corresponds to the motion of the second body
in a central potential of the fixed first body. An additional
integral in this case is $p_{2}$.
\item $m_{2}=\infty$ or equivalently $m=m_{1}$. This case corresponds to the motion of the
first body in a central potential of the fixed second one. In this
case the function $h_{s}$ can be represented in the form
$$
h_{s}=\frac{1}{2m_{1}R^{2}}\left(p_{\theta}+p_{0}\right)^{2}+
\frac{\left(p_{1}\sin\theta+p_{2}\cos\theta\right)^{2}}{2m_{1}R^{2}\sin^{2}\theta}+V(\theta).
$$
Since $\{p_{1}\sin\theta+p_{2}\cos\theta,p_{\theta}+p_{0}\}=0$, an
additional integral in this case is
$p_{1}\sin\theta+p_{2}\cos\theta$.
\end{enumerate}

One can also derive from (\ref{hh3}) the Hamiltonian function for
the restricted two-body problem. Let the first body be a "heavy"
one ($m_{1}\to\infty$) and it moves along a geodesic $\Gamma$ with
a constant velocity $\omega$. Since the momentum $p_{2}$
corresponds to the rotation of the second body around the first
one it holds $p_{2}/m_{1}\to0$. Let $\psi$ be an angle between the
geodesic $\Gamma$ and the geodesic $\Upsilon$ connecting the
bodies. Since the momentum $p_{0}$ corresponds to the motion of
bodies along the geodesic $\Upsilon$ then
$p_{0}/m_{1}\to\omega\cos\psi$. At last the momentum $p_{1}$
corresponds to the motion of the system in the direction,
perpendicular to $\Upsilon$, therefore
$p_{1}/m_{1}\to\omega\sin\psi$. Also one should omit the term
$\gamma^{2}/(2m_{1}R^{2})$, which is independent from time and
tends to infinity, since $\gamma/m_{1}\to\const\neq0$. This
corresponds to the infinite kinetic energy of the "heavy" body as
$m_{1}\to\infty$.

Thus from (\ref{hh3}) one gets the Hamiltonian function for the
restricted two-body problem
$$
h_{s,r}=\frac{1}{2m_{2}R^{2}}\left(p_{\theta}^{2}+\frac{p_{2}^{2}}{\sin^{2}\theta}\right)+
\frac{\omega}{R^{2}}\left(p_{\theta}\cos\psi+p_{2}\sin\psi\cot\theta\right)+V(\theta).
$$
Besides, since
$$\left\{\frac{p_{0}}{m_{1}},p_{2}\right\}\to\omega\{\cos\psi,p_{2}\}=-\omega\sin\psi\{\psi,p_{2}\}\quad
\text{and}\quad
\left\{\frac{p_{0}}{m_{1}},p_{2}\right\}=\frac{p_{1}}{m_{1}}\to\omega\sin\psi$$
one gets in the limiting case $\{\psi,p_{2}\}=-1$. Hence, the
variables $p_{\psi}:=-p_{2},\psi$ are canonical and we obtain
$$
h_{s,r}=\frac{1}{2m_{2}R^{2}}\left(p_{\theta}^{2}+
\frac{p_{\psi}^{2}}{\sin^{2}\theta}\right)+\frac{\omega}{R^{2}}\left(
p_{\theta}\cos\psi-p_{\psi}\sin\psi\cot\theta\right)+V(\theta)
$$
that coincides up to notations with corresponding expressions from
\cite{Zi} and \cite{MacPrz}.

It holds
$$
\left\{p_{\theta}^{2}+
\frac{p_{\psi}^{2}}{\sin^{2}\theta},p_{\theta}\cos\psi-p_{\psi}\sin\psi\cot\theta\right\}=0
$$
that corresponds to the integrability of the free motion.

Note that the Newton and the oscillator potentials mentioned in
the Introduction have the following forms:
\begin{equation}\label{SPot}
V_{N}=-\alpha\cot\theta,\;V_{o}=\frac{\beta}2\tan^{2}\theta,\;\alpha,\beta=\const,\,
\alpha,\beta>0.
\end{equation}

Our main result for the spherical case is the following theorem.
\begin{theore}\label{MainThS}
The complexified Hamiltonian system with Hamiltonian function
(\ref{h1}) and potentials (\ref{SPot}) does not admit an
additional meromorphic first integral in the case $m_{1}\neq m,\,
m_{1}m\alpha\beta\gamma\neq0$, where
$p_{0}^{2}+p_{1}^{2}+p_{2}^{2}=\gamma^{2}$.
\end{theore}

\subsection{The reduced two-body problem on the hyperbolic plane $\mathbf{H}^{2}$}

Let ${\bf H}^{2}$ be the hyperbolic plane with a sectional
curvature $-1/R^{2}$. The configuration space for the two-body
problem on ${\bf H}^{2}$ is $Q=\left({\bf H}^{2}\times{\bf
H}^{2}\right)\backslash\diag$. Here, there are no opposite points
and the phase space $M:=T^{*}Q$ can be represented as
$$
M=T^{*}\mathbb{R}_{+}\times T^{*}\Or_{0}(1,2),
$$
where $\mathbb{R}_{+}:=(0,\infty)$ and $\Or_{0}(1,2)$ is the
identity component of the isometry group for ${\bf H}^{2}$. The
Lie algebra for the Lie group $\Or_{0}(1,2)$ is
$\mathfrak{so}(1,2)$.

The reduced phase space for $M$ is
$$
\widetilde{M}=T^{*}\mathbb{R}_{+}\times\mathcal{O},
$$
where $\mathcal{O}$ is a $\Or_{0}(1,2)$-orbit w.r.t.\ the
coadjoint action in the space $\mathfrak{so}^{*}(1,2)$. It is
endowed with the Kirillov symplectic form.

The Lie algebra $\mathfrak{so}(1,2)$, endowed with Killing form,
is the Minkowski space with the signature $(1,2)$ and can be
naturally identified with its dual space $\mathfrak{so}^{*}(1,2)$.
Let $p_{0},p_{1},p_{2}$ be orthogonal coordinates in
$\mathfrak{so}^{*}(1,2)$ w.r.t.\ the Killing form
$\Kil:=p_{0}^{2}+p_{1}^{2}-p_{2}^{2}$. The orbits
$\mathcal{O}\subset\mathfrak{so}^{*}(1,2)$ are of the following
types
\begin{enumerate}
\item $\mathcal{O}_{0}=(0)$;
\item $\mathcal{O}_{c}=\left((p_{0},p_{1},p_{2})|\,p_{0}^{2}+p_{1}^{2}-p_{2}^{2}=0,
(p_{0},p_{1},p_{2})\ne(0,0,0)\right)$;
\item $\mathcal{O}_{\gamma}=\left((p_{0},p_{1},p_{2})|\,p_{0}^{2}+p_{1}^{2}-p_{2}^{2}=\gamma
\right),\;\gamma>0$;
\item $\mathcal{O}_{\gamma}=\left((p_{0},p_{1},p_{2})|\,p_{0}^{2}+p_{1}^{2}-p_{2}^{2}=\gamma
\right),\;\gamma<0$.
\end{enumerate}
The orbit $\mathcal{O}_{c}$ is a cone without its vertex, the
orbit $\mathcal{O}_{\gamma}$ for $\gamma>0$ is a one-sheet
hyperboloid and for $\gamma<0$ it is a two-sheet hyperboloid.

The two-body Hamiltonian function on $\widetilde{M}$ can be
obtained from (\ref{h1}) by the formal substitution
$$
p_{0}\to ip_{0},\,p_{1}\to ip_{1},\,p_{2}\to p_{2},\,p_{r}\to
ip_{r},\,r\to -ir,\,R\to iR,
$$
where $i$ is the complex unit. This leads to the expression
\begin{equation}\label{hh1}
h_{h}=\frac{(1-r^{2})^{2}}{8mR^{2}}\left(p_{r}^{2}+\frac{p_{2}^{2}}{r^{2}}\right)+
\frac{(1-r^{2})p_{r}p_{0}+\gamma}{2m_{1}R^{2}}+\frac{p_{2}^{2}}{m_{1}R^{2}}+
\frac{(1+r^{2})p_{1}p_{2}}{2m_{1}R^{2}r}+V(r),
\end{equation}
where $r:=\tanh\left(\rho/(2R)\right)$; $\rho\in\mathbb{R}_{+}$ is
the distance between bodies, and
$p_{0}^{2}+p_{1}^{2}-p_{2}^{2}=\gamma$ on the orbit
$\mathcal{O}_{\gamma}$.

The Poisson brackets for variables $r,p_{r},p_{0},p_{1},p_{2}$ are
as follows:
\begin{gather*}
\{r,p_{r}\}=1,\,\{p_{0},p_{1}\}=p_{2},\,\{p_{1},p_{2}\}=-p_{0},\\
\{p_{0},p_{2}\}=p_{1},\,\{r,p_{i}\}=0,\,\{p_{r},p_{i}\}=0,\,i=0,1,2.
\end{gather*}

One can define canonical variables $p_{\varphi},\varphi$ on
$\mathcal{O}_{\gamma}\ne\pt$ by the formulas
$$
p_{0}=\sqrt{\gamma+p_{2}^{2}}\sin\varphi,\,p_{1}=\sqrt{\gamma+p_{2}^{2}}\cos\varphi,
p_{2}=p_{\varphi}.
$$
Then one gets the expression
\begin{align}\label{hh2}
h_{h}&=\frac{(1-r^{2})^{2}}{8mR^{2}}\left(p_{r}^{2}+\frac{p_{\varphi}^{2}}{r^{2}}\right)+
\frac{p_{\varphi}^{2}}{m_{1}R^{2}}
+\frac{\sqrt{\gamma+p_{\varphi}^{2}}}{2m_{1}R^{2}}
\left(p_{\varphi}\frac{1+r^{2}}r\cos\varphi
+(1-r^{2})p_{r}\sin\varphi\right)\nonumber\\&+V(r)+\frac{\gamma}{2m_{1}R^{2}}
\end{align}
that up to notations coincides with formulas (11)-(13) from
\cite{Sh004}.

Let us also define a new pair of canonical variables
$\theta,p_{\theta}$ such that
$$
\theta=\frac{\rho}{R},\,p_{\theta}=\frac12(1-r^{2})p_{r},\,\{\theta,p_{\theta}\}=1.
$$
Then, one gets from (\ref{hh1})
\begin{align}\begin{split}
h_{h}&=\frac{1}{2mR^{2}}\left(p_{\theta}^{2}+\frac{p_{2}^{2}}{\sinh^{2}\theta}\right)+
\frac{p_{\theta}p_{0}}{m_{1}R^{2}}+\frac{p_{2}^{2}}{m_{1}R^{2}}+
\frac{p_{1}p_{2}}{m_{1}R^{2}}\coth\theta+\frac{\gamma}{2m_{1}R^{2}}
+V(\theta).
\end{split}\label{h3}
\end{align}

Cases of a trivial integrability of the reduced two-body problem
on the hyperbolic plane are similar to those described above for
the spherical case and correspond to one of the equalities:
$V(\theta)=0,\,m_{1}=\infty$ or $m=m_{1}$.

The Newton and the oscillator potentials now have the following
forms:
\begin{equation}\label{HPot}
V_{N}=-\alpha\coth\theta,\;V_{o}=\frac{\beta}2\tanh^{2}\theta,\;\alpha,\beta=\const,\,
\alpha,\beta>0.
\end{equation}

The main result of the present paper for the hyperbolic case is as
follows.
\begin{theore}\label{MainThH}
The complexified Hamiltonian system with Hamiltonian function
(\ref{hh1}) and potentials (\ref{HPot}) does not admit an
additional meromorphic first integral in the case $m_{1}\neq m,\,
m_{1}m\alpha\beta\gamma\neq0,\,\gamma>0$, where
$p_{0}^{2}+p_{1}^{2}-p_{2}^{2}=\gamma$.
\end{theore}

\section{The result of the Morales-Ramis theory}\label{M-RT}\markright{\ref{M-RT}
The result of the Morales-Ramis theory}

Here we present the result from the Morales-Ramis theory
\cite{MorRu}, which will be used for proving the absence of an
additional meromorphic integral for the reduced two-body problem
on ${\bf S}^{2}$ and ${\bf H}^{2}$ with potentials (\ref{SPot})
and (\ref{HPot}).

Let $M$ be a complex analytic manifold ($\dim_{\mathbb{C}}M=n$)
and
\begin{equation}\label{DES}
\frac{dx}{dt}=v(x),\,t\in\mathbb{C},x\in M
\end{equation}
be a system of differential equations, where $v$ is a holomorphic
vector field on $M$. Let $x=\varphi(t)$ be a particular
nonconstant solution of (\ref{DES}) and $\Gamma$ be the Riemann
surface in $M$, defined by the maximal analytic continuation of
$\varphi(t)$. The {\it variational equations} along $\varphi(t)$
have the form
\begin{equation}\label{VE}
\frac{d\xi}{dt}=V(v)\xi,\;\xi\in T_{\Gamma}M,
\end{equation}
where $T_{\Gamma}M$ is the restriction of the tangent bundle $TM$
onto $\Gamma$ and $$V(v):=\frac{\partial v}{\partial
x}\left(\varphi(t)\right)$$ is a morphism $T_{\Gamma}M\to
T_{\Gamma}M$ of the vector bundle $T_{\Gamma}M$.

Let system (\ref{DES}) be a Hamiltonian one with a Hamiltonian
function $H$; in particular, $n$ is even. Then the order of
(\ref{VE}) can be reduced by two. Indeed, let $M_{\varepsilon}$ be
a submanifold in $M$ defined as $M_{\varepsilon}:=\left(x\in
M|\,H(x)=\varepsilon\right),\,\varepsilon=\const$;
$N:=T_{\Gamma}M_{\varepsilon}/T\Gamma$ be the normal bundle of the
surface $\Gamma$ in $TM_{\varepsilon}$, and
$\pi:\;T_{\Gamma}M_{\varepsilon}\to N$ be the canonical
projection. Note that $\dim_{\mathbb{C}}N=n-2$. Then
\begin{equation}\label{NVE}
\frac{d\eta}{dt}=\pi\left[V(v)\left(\pi^{-1}\eta\right)\right],\;\eta\in
N
\end{equation}
is a well-defined system of differential equations on $N$ since
$V(v)\left(T\Gamma\right)\subset T\Gamma$ and $\pi(T\Gamma)$ is
the null section of $N$. System (\ref{NVE}) is called the {\it
normal variational equations}.

Let $\EuScript{G}$ be the differential Galois group \cite{Kap},
\cite{vdPS} for system (\ref{NVE}), i.e.\ a matrix group acting on
fundamental solutions of (\ref{NVE}) which does not change
polynomial relations between them. Let also $\EuScript{G}_{0}$ be
the identity component for $\EuScript{G}$.

\begin{theore}[see \cite{MorRu}, theorem 4.1]\label{MRTh}
Suppose that there are $n/2$ meromorphic first integrals of
Hamiltonian system (\ref{DES}) that are in involution and are
independent in some neighborhood of $\Gamma$. Then
$\EuScript{G}_{0}$ is an abelian group.
\end{theore}

\section{Particular solutions and variational equations}\label{ParSol}
\markright{\ref{ParSol} Particular solutions and variational
equations}

In order to simplify notations one can multiply Hamiltonian
functions (\ref{hh3}) and (\ref{h3}) by $m_{1}R^{2}$ that is
equivalent to changing the scale of the time axis and omit the
constant summand. Thus one gets
\begin{gather}\label{hsSimple}
h_{s}=\frac{1}{2\mu}\left(p_{\theta}^{2}+\frac{p_{2}^{2}}{\sin^{2}\theta}\right)+
p_{\theta}p_{0}-p_{2}^{2}+p_{1}p_{2}\cot\theta+V(\theta),\\
\label{hhSimple}
h_{h}=\frac{1}{2\mu}\left(p_{\theta}^{2}+\frac{p_{2}^{2}}{\sinh^{2}\theta}\right)+
p_{\theta}p_{0}+p_{2}^{2}+p_{1}p_{2}\coth\theta+V(\theta),
\end{gather}
where $\mu:=m/m_{2}=m_{1}/(m_{1}+m_{2})\neq0$.

Consider Hamiltonian systems with Hamiltonian functions
(\ref{hsSimple}) and (\ref{hhSimple}) on reduced manifolds
$\widetilde{M}_{\gamma}=T^{*}I\times\mathcal{O}_{\gamma}$ in the
spherical case and
$\widetilde{M}_{\gamma}=T^{*}\mathbb{R}_{+}\times\mathcal{O}_{\gamma}$
in the hyperbolic case. For any potential $V(\theta)$ there are
trajectories defined by the equalities
$p_{0}=p=\const\neq0,\,p_{1}=p_{2}=0$. They correspond to the
bodies motion along a common geodesic. For the spherical case all
nondegenerate manifolds $\widetilde{M}_{\gamma},\,\gamma>0$
contain such a trajectory and for the hyperbolic case only those
with $\gamma=p^{2}>0$. Denote the maximal analytic continuation of
this trajectory by $\Gamma$ in accordance with section \ref{M-RT}.

One can choose $p_{1}$ and $p_{2}$ as local coordinates in a
neighborhood of $\Gamma$. Then using (\ref{PoissonS}) and
(\ref{evolution}) one gets the normal variational equations in the
spherical case:
\begin{align}\begin{split}
\frac{dp_{1}}{dt}&=-p\cot\theta
p_{1}+\left(2p+p_{\theta}-\frac{p}{\mu\sin^{2}\theta}\right)p_{2},\\
\frac{dp_{2}}{dt}&=-p_{\theta}p_{1}+p\cot\theta p_{2},
\end{split}\label{NVES}
\end{align}
where $p_{\theta}=p_{\theta}(t),\,\theta=\theta(t)$ is a solution
of the Hamiltonian system with the Hamiltonian function
\begin{equation}\label{GammaHamF}
h_{0}=\frac{1}{2\mu}p_{\theta}^{2}+
pp_{\theta}+V(\theta)=\frac{1}{2\mu}\left(p_{\theta}+\mu
p\right)^{2}+V(\theta)-\frac{\mu}2p^{2}.
\end{equation}
The normal variational equations in the hyperbolic case are
\begin{align}\begin{split}
\frac{dp_{1}}{dt}&=-p\coth\theta
p_{1}-\left(2p+p_{\theta}+\frac{p}{\mu\sinh^{2}\theta}\right)p_{2},\\
\frac{dp_{2}}{dt}&=-p_{\theta}p_{1}+p\coth\theta p_{2},
\end{split}\label{NVEH}
\end{align}
where again $p_{\theta}=p_{\theta}(t),\,\theta=\theta(t)$ is a
solution of the Hamiltonian system with the Hamiltonian function
(\ref{GammaHamF}).

One can compare these normal variational systems with their
analogues for the restricted two-body problem from \cite{Zi} and
\cite{MacPrz}. For example, in the spherical case latter system
can be written as
\begin{align}\begin{split}
\frac{dp_{1}}{dt}&=-\omega\cot\theta
p_{1}+\frac{p_{2}}{\sin^{2}\theta},\\ \frac{dp_{2}}{dt}&=\omega
p_{\theta}p_{1}+\omega\cot\theta p_{2}.
\end{split}\label{NVER} \end{align}

The key factor for determination of a differential Galois group
for a system of linear differential equations is its reducibility
to a system with rational coefficients. For the Newton and
oscillator potentials, such reduction for systems (\ref{NVES}),
(\ref{NVEH}) and (\ref{NVER}) is possible.

\subsection{The Newton potential}\label{ParSol1}

For the Newton potential $V=V_{N}=-\alpha\cot\theta$, system
(\ref{NVER}) becomes a Fuchsian one w.r.t.\ the independent
variable $p_{\theta}$ that was found in \cite{Zi}. This fact is
also valid for system (\ref{NVES}).

Indeed, denote $z:=(p_{\theta}+\mu p)/\alpha$. One can easily
check that system (\ref{NVES}) can be written as
\begin{align}\begin{split}
p_{1}'(z)&=A(z)p_{1}+B(z)p_{2}\\ p_{2}'(z)&=C(z)p_{1}-A(z)p_{2}
\end{split}\label{NVESz}\end{align}
with respect to the independent variable $z$, where
\begin{align*}
A(z)=\frac{pf(z)}{1+f^{2}(z)},\;B(z)=\frac{p}{\mu}-\frac{\alpha
z+(2-\mu)p}{1+f^{2}(z)},\;C(z)=\frac{\alpha z-\mu
p}{1+f^{2}(z)},\; f(z)=\frac{\alpha z^{2}}{2\mu}-\varepsilon
\end{align*}
and the trajectory $\Gamma$ corresponds to the equation
$$
\frac{\alpha^{2}z^{2}}{2\mu}-\alpha\cot\theta=\alpha\varepsilon=\const.
$$
System (\ref{NVESz}) is Fuchsian (see appendix) with five regular
singular points
$z_{1,2}=\pm\varkappa,\,z_{3,4}=\pm\lambda,\,z_{5}=\infty$, where
$$
\varkappa:=\sqrt{\frac{2\mu}{\alpha}(\varepsilon+i)}\notin\mathbb{R},\;
\lambda:=\sqrt{\frac{2\mu}{\alpha}(\varepsilon-i)}\notin\mathbb{R},\quad
\text{for}\quad \varepsilon\in\mathbb{R}.
$$
We shall express all coefficients through four parameters
$p,\mu,\varkappa,\lambda$. In particular it holds
$$
f(z)=i\frac{2z^{2}-\lambda^{2}-\varkappa^{2}}{\varkappa^{2}-\lambda^{2}},\;
f^{2}(z)+1=2i\frac{(z^{2}-\lambda^{2})(z^{2}-\varkappa^{2})}{\varkappa^{2}-\lambda^{2}}.
$$
One can transform (\ref{NVESz}) into the linear differential
equation for $p_{2}(z)$ of the second order
\begin{equation}\label{p21}
p_{2}''(z)=\frac{C'}Cp'_{2}+\left(\frac{C'}CA+A^{2}+CB-A'\right)p_{2}
\end{equation}
and then into equation (\ref{Diff2Equation2}) for the function
$y(z):=p_{2}(z)/\sqrt{C}$, where
\begin{equation*}
r(z)=\frac{C'}CA+A^{2}+CB-A'-\frac12\left(\frac{C'}C\right)'+\frac14\left(\frac{C'}C\right)^{2}.
\end{equation*}
For evaluation of the function $r(z)$ one can use computer
analytical calculations, which lead to
\begin{equation}\label{r(z)}
r(z)=\sum_{j=1}^{4}\left(\frac{\alpha_{j}}{(z-z_{j})^{2}}+
\frac{\beta_{j}}{z-z_{j}}\right)+\frac3{4(z-z_{0})^{2}}=
\frac3{4z^{2}}+O(\frac1{z^{3}})\quad\text{as}\quad z\to\infty,
\end{equation}
where
\begin{align}\label{r(z)SN}
\alpha_{1}&=\frac{1-\mu}{64\varkappa^{2}}\left(p(\mu-1)(\lambda^{2}-\varkappa^{2})+
4i\varkappa(\mu+1)\right)\left(p(\lambda^{2}-\varkappa^{2})+4i\varkappa\right),\notag\\
\alpha_{2}&=\frac{1-\mu}{64\varkappa^{2}}\left(p(\mu-1)(\lambda^{2}-\varkappa^{2})-
4i\varkappa(\mu+1)\right)\left(p(\lambda^{2}-\varkappa^{2})-4i\varkappa\right),\notag\\
\alpha_{3}&=\frac{1-\mu}{64\lambda^{2}}\left(p(\mu-1)(\varkappa^{2}-\lambda^{2})-
4i\lambda(\mu+1)\right)\left(p(\varkappa^{2}-\lambda^{2})-4i\lambda\right),\notag\\
\alpha_{4}&=\frac{1-\mu}{64\lambda^{2}}\left(p(\mu-1)(\varkappa^{2}-\lambda^{2})+
4i\lambda(\mu+1)\right)\left(p(\varkappa^{2}-\lambda^{2})+4i\lambda\right),\notag\\
\beta_{1}&=\frac{\mu-1}{64(\varkappa^{2}-\lambda^{2})\varkappa^{3}}
\left((\mu-1)(5\varkappa^{2}-\lambda^{2})(\varkappa^{2}-\lambda^{2})^{2}p^{2}-32i\mu
(\varkappa^{2}-\lambda^{2})\varkappa^{3}p\right.\notag\\&\left.-16(\mu+1)\varkappa^{2}(3\varkappa^{2}+\lambda^{2})
\right),\\
\beta_{2}&=\frac{\mu-1}{64(\varkappa^{2}-\lambda^{2})\varkappa^{3}}
\left((1-\mu)(5\varkappa^{2}-\lambda^{2})(\varkappa^{2}-\lambda^{2})^{2}p^{2}-32i\mu
(\varkappa^{2}-\lambda^{2})\varkappa^{3}p\right.\notag\\&\left.+16(\mu+1)\varkappa^{2}(3\varkappa^{2}+\lambda^{2})
\right),\notag\\
\beta_{3}&=\frac{\mu-1}{64(\varkappa^{2}-\lambda^{2})\lambda^{3}}
\left((1-\mu)(5\lambda^{2}-\varkappa^{2})(\varkappa^{2}-\lambda^{2})^{2}p^{2}+32i\mu
(\varkappa^{2}-\lambda^{2})\lambda^{3}p\right.\notag\\&\left.+16(\mu+1)\lambda^{2}(3\lambda^{2}+\varkappa^{2})
\right),\notag\\
\beta_{4}&=\frac{\mu-1}{64(\varkappa^{2}-\lambda^{2})\lambda^{3}}
\left((\mu-1)(5\lambda^{2}-\varkappa^{2})(\varkappa^{2}-\lambda^{2})^{2}p^{2}+32i\mu
(\varkappa^{2}-\lambda^{2})\lambda^{3}p\right.\notag\\&\left.-16(\mu+1)\lambda^{2}(3\lambda^{2}+\varkappa^{2})
\right),\; z_{0}=\frac{\mu
p}{\alpha}=\frac{p(\varkappa^{2}-\lambda^{2})}{4i}.\notag
\end{align}
For $\mu=1$ the expression for $r(z)$ is very simple
\begin{equation}\label{rForMu=1}
r(z)=\frac3{4(z-z_{0})^{2}}.
\end{equation}

\begin{Lem}\label{L1}
Suppose that
$\alpha,\varepsilon,\mu,p\in\mathbb{R},\,\mu\neq0,1,\,p\neq0$ and
\begin{equation}\label{NonEq}
(\sqrt{\varepsilon^{2}+1}-\varepsilon)(\varepsilon^{2}+1)\neq
\frac{(\mu-1)^{2}p^{2}}{4\alpha\mu};
\end{equation}
then $\alpha_{i}\notin\mathbb{R}$ for $i=1,2,3,4$.
\end{Lem}
\begin{proof}
Direct calculations imply
$\alpha_{1}=\dfrac{\mu^{2}-1}4-\dfrac{\mu-1}{4\alpha}p\mu^{2}\left(\dfrac2\varkappa+
\dfrac{p(1-\mu)}{\varkappa^{2}\alpha}\right)$ and
$$
\frac2\varkappa+\frac{p(1-\mu)}{\varkappa^{2}\alpha}=\frac{\sqrt{2\alpha}}
{\sqrt{\mu}\sqrt{\varepsilon^{2}+1}}\left(\sqrt{\varepsilon-i}+
\frac{(\mu-1)pi}{2\sqrt{2\alpha\mu}\sqrt{\varepsilon^{2}+1}}\right)+
\frac{(1-\mu)p\varepsilon}{2\mu(\varepsilon^{2}+1)}.
$$
Therefore $\alpha_{1}\notin\mathbb{R}$ iff
$$
\im\sqrt{\varepsilon-i}=\pm\frac1{\sqrt{2}}\sqrt{\sqrt{\varepsilon^{2}+1}-\varepsilon}\neq
\frac{(1-\mu)p}{2\sqrt{2\alpha\mu}\sqrt{\varepsilon^{2}+1}}
$$
that is equivalent to (\ref{NonEq}). The consideration for
$\alpha_{2},\alpha_{3}$ and $\alpha_{4}$ is similar.
\end{proof}

In the hyperbolic case, system (\ref{NVEH}) again is reduced to
the Fuchsian system (\ref{NVESz}), where now
\begin{align*}
A(z)=\frac{pf(z)}{f^{2}(z)-1},\;B(z)=\frac{p}{\mu}+\frac{\alpha
z+(2-\mu)p}{f^{2}(z)-1},\;C(z)=\frac{\alpha z-\mu
p}{f^{2}(z)-1},\; f(z)=\frac{\alpha z^{2}}{2\mu}-\varepsilon.
\end{align*}
The trajectory $\Gamma$ corresponds here to the equation
$$
\frac{\alpha z^{2}}{2\mu}-\coth\theta=\varepsilon.
$$
In this case singular points are
$z_{1,2}=\pm\varkappa,\,z_{3,4}=\pm\lambda,z_{5}=\infty$ for
$$
\varkappa:=\sqrt{\frac{2\mu}{\alpha}(\varepsilon+1)},\;
\lambda:=\sqrt{\frac{2\mu}{\alpha}(\varepsilon-1)}.
$$
Again for the function $y(z)=p_{2}(z)/\sqrt{C(z)}$ one gets
equation (\ref{Diff2Equation2}) with $r(z)$ given by (\ref{r(z)}),
where
\begin{align}\label{r(z)HN}
\alpha_{1}&=\frac{\mu-1}{64\varkappa^{2}}\left(p(\mu-1)(\varkappa^{2}-\lambda^{2})-
4\varkappa(\mu+1)\right)\left(p(\varkappa^{2}-\lambda^{2})-4\varkappa\right),\notag\\
\alpha_{2}&=\frac{\mu-1}{64\varkappa^{2}}\left(p(\mu-1)(\varkappa^{2}-\lambda^{2})+
4\varkappa(\mu+1)\right)\left(p(\varkappa^{2}-\lambda^{2})+4\varkappa\right),\notag\\
\alpha_{3}&=\frac{\mu-1}{64\lambda^{2}}\left(p(\mu-1)(\varkappa^{2}-\lambda^{2})-
4\lambda(\mu+1)\right)\left(p(\varkappa^{2}-\lambda^{2})-4\lambda\right),\notag\\
\alpha_{4}&=\frac{\mu-1}{64\lambda^{2}}\left(p(\mu-1)(\varkappa^{2}-\lambda^{2})+
4\lambda(\mu+1)\right)\left(p(\varkappa^{2}-\lambda^{2})+4\lambda\right),\notag\\
\beta_{1}&=\frac{\mu-1}{64(\varkappa^{2}-\lambda^{2})\varkappa^{3}}
\left((\mu-1)(\lambda^{2}-5\varkappa^{2})(\varkappa^{2}-\lambda^{2})^{2}p^{2}+32\mu
(\varkappa^{2}-\lambda^{2})\varkappa^{3}p\right.\notag\\&\left.-16(\mu+1)\varkappa^{2}
(3\varkappa^{2}+\lambda^{2})\right),\\
\beta_{2}&=\frac{\mu-1}{64(\varkappa^{2}-\lambda^{2})\varkappa^{3}}
\left((\mu-1)(5\varkappa^{2}-\lambda^{2})(\varkappa^{2}-\lambda^{2})^{2}p^{2}+32\mu
(\varkappa^{2}-\lambda^{2})\varkappa^{3}p\right.\notag\\&\left.+16(\mu+1)\varkappa^{2}(3\varkappa^{2}+
\lambda^{2})\right),\notag\\
\beta_{3}&=\frac{\mu-1}{64(\varkappa^{2}-\lambda^{2})\lambda^{3}}
\left((\mu-1)(5\lambda^{2}-\varkappa^{2})(\varkappa^{2}-\lambda^{2})^{2}p^{2}-32\mu
(\varkappa^{2}-\lambda^{2})\lambda^{3}p\right.\notag\\&\left.+16(\mu+1)\lambda^{2}
(3\lambda^{2}+\varkappa^{2})\right),\notag\\
\beta_{4}&=\frac{\mu-1}{64(\varkappa^{2}-\lambda^{2})\lambda^{3}}
\left((\mu-1)(\varkappa^{2}-5\lambda^{2})(\varkappa^{2}-\lambda^{2})^{2}p^{2}-32\mu
(\varkappa^{2}-\lambda^{2})\lambda^{3}p\right.\notag\\&\left.-16(\mu+1)\lambda^{2}(3\lambda^{2}+
\varkappa^{2})\right),\; z_{0}=\frac{\mu
p}{\alpha}=\frac{p(\varkappa^{2}-\lambda^{2})}{4}.\notag
\end{align}
For $\mu=1$ the expression for $r(z)$ coincides with
(\ref{rForMu=1}).
\begin{Lem}\label{L2}
Suppose that
$\alpha,\varepsilon,\mu,p\in\mathbb{R},\,\mu\neq0,1;\,p\neq0,\,\varepsilon<-1$;
then $\alpha_{i}\notin\mathbb{R}$ for $i=1,2,3,4$.
\end{Lem}
\begin{proof}
Evidently, $\varkappa,\lambda\in i\mathbb{R}\backslash(0)$ and one
gets
$$
i\im\alpha_{1}=-i\im\alpha_{2}=\frac{\mu^{2}(1-\mu)p}{2\alpha\varkappa}\neq0,\quad
i\im\alpha_{3}=-i\im\alpha_{4}=\frac{\mu^{2}(1-\mu)p}{2\alpha\lambda}\neq0.
$$
\end{proof}

\subsection{The oscillator potential}\label{ParSol2}

As above, using the independent variable $z:=(p_{\theta}+\mu
p)/\beta$ for the oscillator potential $V=\beta\tan^{2}\theta/2$
one can reduce system (\ref{NVES}) to the system
\begin{align}\begin{split}
p_{1}'(z)&=A(z)p_{1}+B(z)\sqrt{f(z)}p_{2}\\
p_{2}'(z)&=C(z)\sqrt{f(z)}p_{1}-A(z)p_{2}
\end{split}\label{NVESzOs}\end{align}
with coefficients
\begin{align*}
A(z)=\frac{p}{f(f+1)},\;B(z)=\frac{p}{\mu f^{2}}-\frac{\beta
z+(2-\mu)p}{f(f+1)},\;C(z)=\frac{\beta z-\mu p}{f(f+1)},
\end{align*}
where $f(z)=\tan^{2}\theta=-\dfrac{\beta}{\mu}
z^{2}+2\varepsilon$.

In the general case coefficients of system (\ref{NVESzOs}) are not
rational due to the appearance of $\sqrt{f(z)}$. The same
difficulty for the restricted two-body problem was overcome in
\cite{Zi2} by assumption $\varepsilon=0$, when
$\sqrt{f(z)}=\sqrt{-\beta/\mu}z$. On the other hand it was noted
in \cite{MacPrz} that one can pass on to a second order
differential equation with rational coefficients.

Using the latter approach one gets from (\ref{NVESzOs}) the
following equation for $p_{2}(z)$:
\begin{equation*}\label{p22}
p_{2}''(z)=\left(\frac{C'}C+\frac{f'}{2f}\right)p'_{2}+
\left(\left(\frac{C'}C+\frac{f'}{2f}\right)A+A^{2}+CBf-A'\right)p_{2},
\end{equation*}
which can be reduced to equation (\ref{Diff2Equation2}) by the
substitution $p_{2}(z)=y(z)\sqrt{C(z)}\left(f(z)\right)^{1/4}$.
Here
\begin{equation*}
r(z)=\left(\frac{C'}C+\frac{f'}{2f}\right)A+A^{2}+CBf-A'-\frac12\left(\frac{C'}C+
\frac{f'}{2f}\right)'+\frac14\left(\frac{C'}C+\frac{f'}{2f}\right)^{2}.
\end{equation*}

Denote by
$z_{0}:=p\mu/\beta,\,z_{1,2}:=\pm\varkappa:=\pm\sqrt{\mu(2\varepsilon+1)/\beta}$
and $z_{3,4}:=\pm\lambda:=\pm\sqrt{2\mu\varepsilon/\beta}$ zeros
of functions $C(z),\,f(z)+1$ and $f(z)$ respectively. Using
computer calculations one gets
\begin{equation}\label{r(z)Osc}
r(z)=\sum_{j=0}^{4}\left(\frac{\alpha_{j}}{(z-z_{j})^{2}}+
\frac{\beta_{j}}{z-z_{j}}\right)=O(\frac1{z^{4}})\quad\text{as}\quad
z\to\infty,
\end{equation}
where
\begin{align}\label{coeffr(z)SO}
\alpha_{0}&=\frac34,\;\alpha_{3}=\alpha_{4}=-\frac3{16},\;
\beta_{0}=\frac{3(\lambda^{2}-\varkappa^{2})p}{2\left((\lambda^{2}-\varkappa^{2})^{2}p^{2}-
\lambda^{2}\right)},\notag\\
\alpha_{1}&=\frac{\mu-1}{4\varkappa^{2}}\left(p(\mu-1)(\varkappa^{2}-\lambda^{2})-
\varkappa(\mu+1)\right)\left(p(\varkappa^{2}-\lambda^{2})-\varkappa\right),\notag\\
\alpha_{2}&=\frac{\mu-1}{4\varkappa^{2}}\left(p(\mu-1)(\varkappa^{2}-\lambda^{2})+
\varkappa(\mu+1)\right)\left(p(\varkappa^{2}-\lambda^{2})+\varkappa\right),\notag\\
\beta_{1}&=\frac{\mu-1}{4(\varkappa^{2}-\lambda^{2})\varkappa^{3}}
\left((\mu-1)(\lambda^{2}-3\varkappa^{2})(\varkappa^{2}-\lambda^{2})^{2}p^{2}+4\mu
(\varkappa^{2}-\lambda^{2})\varkappa^{3}p\right.\notag\\&\left.-(\mu+1)\varkappa^{2}
(\varkappa^{2}+\lambda^{2})\right),\\
\beta_{2}&=\frac{\mu-1}{4(\varkappa^{2}-\lambda^{2})\varkappa^{3}}
\left((\mu-1)(3\varkappa^{2}-\lambda^{2})(\varkappa^{2}-\lambda^{2})^{2}p^{2}+4\mu
(\varkappa^{2}-\lambda^{2})\varkappa^{3}p\right.\notag\\&\left.+(\mu+1)\varkappa^{2}
(\varkappa^{2}+\lambda^{2})\right),\notag\\
\beta_{3}&=\left[8(\mu-1)^{2}(\varkappa^{2}-\lambda^{2})^{3}p^{3}-
8(3\mu-1)(\mu-1)\lambda(\varkappa^{2}-\lambda^{2})^{2}p^{2}
+\lambda((5-8\mu^{2})\lambda^{2}+3\varkappa^{2})\right.\notag\\&\left.+
(\varkappa^{2}-\lambda^{2})(9\varkappa^{2}+(24\mu^{2}-16\mu-17)\lambda^{2})p
\right]/\left[16\lambda(\varkappa^{2}-\lambda^{2})
((\varkappa^{2}-\lambda^{2})p-\lambda)\right],\notag\\
\beta_{4}&=\left[8(\mu-1)^{2}(\lambda^{2}-\varkappa^{2})^{3}p^{3}-
8(3\mu-1)(\mu-1)\lambda(\varkappa^{2}-\lambda^{2})^{2}p^{2}
+\lambda((5-8\mu^{2})\lambda^{2}+3\varkappa^{2})\right.\notag\\&\left.+
(\lambda^{2}-\varkappa^{2})\left(9\varkappa^{2}+(24\mu^{2}-16\mu-17)\lambda^{2}\right)p
\right]/\left[16\lambda(\varkappa^{2}-\lambda^{2})
((\varkappa^{2}-\lambda^{2})p+\lambda)\right].\notag
\end{align}
For $\mu=1$ the function $r(z)$ has the form
\begin{equation}\label{rOscMu=1}
r(z)=\frac34\frac{(z_{0}^{2}-2\lambda^{2})z^{2}+2\lambda^{2}z_{0}z+
\lambda^{2}(\lambda^{2}-2z_{0}^{2})}{(z-z_{0})^{2}(z^{2}-\lambda^{2})^{2}}.
\end{equation}

In the hyperbolic case system (\ref{NVEH}) again is reduced to the
system (\ref{NVESzOs}), where now
\begin{align*}
A(z)=\frac{p}{f(1-f)},\;B(z)=\frac{p}{\mu f^{2}}+\frac{\beta
z+(2-\mu)p}{f(1-f)},\;C(z)=\frac{\beta z-\mu p}{f(1-f)},
\end{align*}
and $f(z)=\tanh^{2}\theta=-\dfrac{\beta}{\mu} z^{2}+2\varepsilon$.

Reasoning as above for the spherical case one can get equation
(\ref{Diff2Equation2}) for the function
$y(z):=p_{2}(z)\left(C(z)\right)^{-1/2}\left(f(z)\right)^{-1/4}$.
Now the singular points are
$z_{0}=p\mu/\beta,\,z_{1,2}=\pm\varkappa,\,z_{3,4}=\pm\lambda,z_{5}=\infty$
for
$$
\varkappa:=\sqrt{\frac{\mu}{\beta}(2\varepsilon-1)},\;
\lambda:=\sqrt{\frac{2\mu\varepsilon}{\beta}}
$$
and $r(z)$ is given by (\ref{r(z)Osc}) with
\begin{align}\label{r(z)HO}
\alpha_{0}&=\frac34,\;\alpha_{3}=\alpha_{4}=-\frac3{16},\;
\beta_{0}=\frac{3(\varkappa^{2}-\lambda^{2})p}{2\left((\lambda^{2}-\varkappa^{2})^{2}p^{2}-
\lambda^{2}\right)},\notag\\
\alpha_{1}&=\frac{\mu-1}{4\varkappa^{2}}\left(p(\mu-1)(\varkappa^{2}-\lambda^{2})+
\varkappa(\mu+1)\right)\left(p(\varkappa^{2}-\lambda^{2})+\varkappa\right),\notag\\
\alpha_{2}&=\frac{\mu-1}{4\varkappa^{2}}\left(p(\mu-1)(\varkappa^{2}-\lambda^{2})-
\varkappa(\mu+1)\right)\left(p(\varkappa^{2}-\lambda^{2})-\varkappa\right),\notag\\
\beta_{1}&=\frac{\mu-1}{4(\varkappa^{2}-\lambda^{2})\varkappa^{3}}
\left((\mu-1)(\lambda^{2}-3\varkappa^{2})(\varkappa^{2}-\lambda^{2})^{2}p^{2}-4\mu
(\varkappa^{2}-\lambda^{2})\varkappa^{3}p\right.\notag\\&\left.-(\mu+1)\varkappa^{2}
(\varkappa^{2}+\lambda^{2})\right),\\
\beta_{2}&=\frac{\mu-1}{4(\varkappa^{2}-\lambda^{2})\varkappa^{3}}
\left((\mu-1)(3\varkappa^{2}-\lambda^{2})(\varkappa^{2}-\lambda^{2})^{2}p^{2}-4\mu
(\varkappa^{2}-\lambda^{2})\varkappa^{3}p\right.\notag\\&\left.+(\mu+1)\varkappa^{2}
(\varkappa^{2}+\lambda^{2})\right),\notag\\
\beta_{3}&=\left[8(\mu-1)^{2}(\varkappa^{2}-\lambda^{2})^{3}p^{3}+
8(3\mu-1)(\mu-1)\lambda(\varkappa^{2}-\lambda^{2})^{2}p^{2}
+\lambda((8\mu^{2}-5)\lambda^{2}-3\varkappa^{2})\right.\notag\\&\left.+
(\varkappa^{2}-\lambda^{2})\left(9\varkappa^{2}+(24\mu^{2}-16\mu-17)\lambda^{2}\right)p
\right]/\left[16\lambda(\varkappa^{2}-\lambda^{2})
((\varkappa^{2}-\lambda^{2})p+\lambda)\right],\notag\\
\beta_{4}&=\left[8(\mu-1)^{2}(\lambda^{2}-\varkappa^{2})^{3}p^{3}+
8(3\mu-1)(\mu-1)\lambda(\varkappa^{2}-\lambda^{2})^{2}p^{2}
+\lambda((8\mu^{2}-5)\lambda^{2}-3\varkappa^{2})\right.\notag\\&\left.+
(\lambda^{2}-\varkappa^{2})(9\varkappa^{2}+(24\mu^{2}-16\mu-17)\lambda^{2})p
\right]/\left[16\lambda(\varkappa^{2}-\lambda^{2})
((\varkappa^{2}-\lambda^{2})p-\lambda)\right].\notag
\end{align}
For $\mu=1$ the function $r(z)$ coincides with (\ref{rOscMu=1}).

\begin{Lem}\label{L3}
Suppose that
$\alpha,\varepsilon,\mu,p\in\mathbb{R},\,p\neq0,\,\mu\neq0,1$, and
$\varepsilon<-1/2$ in the spherical case or $\varepsilon<1/2$ in
the hyperbolic case; then $\alpha_{i}\notin\mathbb{R}$ for
$i=1,2$.
\end{Lem}
\begin{proof}
Clearly, in both cases $\varkappa\in i\mathbb{R}\backslash(0)$ and
one gets
$$
i\im\alpha_{1}=-i\im\alpha_{2}=\frac{\mu^{2}(1-\mu)p}{2\beta\varkappa}\neq0.
$$
\end{proof}

\section{Proof of nonintegrability}\label{PNI}\markright{\ref{PNI}
Proof of nonintegrability}

\begin{Lem}\label{L5}
\begin{enumerate}
\item Suppose that assumptions of lemma \ref{L1} are valid. Then the
identity component $\EuScript{G}_{0}$ of the Galois group for
equation (\ref{Diff2Equation2}) with $r(z)$ given by (\ref{r(z)})
and (\ref{r(z)SN}) is not Abelian.
\item Suppose that assumptions of lemma \ref{L2} are valid. Then the
identity component $\EuScript{G}_{0}$ of the Galois group for
equation (\ref{Diff2Equation2}) with $r(z)$ given by (\ref{r(z)})
and (\ref{r(z)HN}) is not Abelian.
\end{enumerate} \end{Lem}
\begin{proof}
We shall prove both claims of this lemma simultaneously. Here,
there are six regular singular points of order $2$:
$z_{j},\,j=0,\ldots,4$ and $z_{5}=\infty$. The difference of
exponents at points $z_{j},\,j=0,\ldots,5$ are
$\Delta_{0}=\Delta_{\infty}=2,\,\Delta_{j}=\sqrt{1+4\alpha_{j}},\,j=1,2,3,4$
and due to lemma \ref{L1} it holds
$\Delta_{j}\not\in\mathbb{R},\,j=1,2,3,4$. Therefore the third
case from lemma \ref{Lem1} is impossible.

Consider the first case of lemma \ref{Lem1}. Here, one or two
linear independent solutions $y_{k}(z)$ of (\ref{Diff2Equation2})
should exist such that $y_{k}'/y_{k}\in\mathbb{C}(z)$.

The rational function $y'_{k}/y_{k}$ has no poles of order more
than one since at such a pole the growth of $y_{k}$ is exponential
that is impossible in the Fuchsian case. Due to the same reason it
should be $\left(y'_{k}/y_{k}\right)(z)\to0$ as $z\to\infty$. This
yields
$$
\frac{y'_{k}(z)}{y_{k}(z)}=\sum_{l}\frac{\delta_{l}}{z-\tilde
z_{l}},\;\delta_{l}\in\mathbb{C}
$$
and one can conclude that
\begin{equation}\label{F}
y_{k}(z)=P_{k}(z)\prod_{j=0}^{4}(z-z_{j})^{\rho_{k}^{(j)}},\,P_{k}(z)\in\mathbb{C}[z],
\end{equation}
where $$\rho_{k}^{(0)}\in\left(-\frac12,\frac32\right),\,
\rho_{k}^{(j)}\in\left(\frac12(1+\Delta_{j}),\frac12(1-\Delta_{j})\right),
\,j=1,\ldots,4,\,k=1,2.$$

Suppose first that there are two such linear independent solutions
$y_{1}$ and $y_{2}$. Then due to lemma \ref{Lem2} it holds
$v(z):=y_{1}(z)y_{2}(z)\in\mathbb{C}(z)$.

Hence, possible exponents for $v(z)$ are $-1,1,3$ at the point
$z_{0}$ and $1,1\pm\sqrt{1+4\alpha_{j}}\not\in\mathbb{R}$ at
points $z_{j},\,j=1,\ldots,4$. The inclusion
$v(z)\in\mathbb{C}(z)$ implies therefore only two possibilities:
$$
v(z)=\frac{P(z)}{z-z_{0}}\prod_{j=1}^{4}(z-z_{j})\quad\text{or}\quad
v(z)=P(z)\prod_{j=1}^{4}(z-z_{j}),\;P(z)\in\mathbb{C}[z].
$$
From (\ref{r(z)}) one can find that exponents for
(\ref{Diff2Equation2}) at $\infty$ are $-3/2$ and $1/2$.
Consequently, the function $v(z)$ grows as $z\to\infty$ no faster
than $z^{3}$. Thus the only possibility (up to a constant nonzero
multiple) for $v$ is
\begin{equation}\label{v(z)}
v(z)=\frac1{z-z_{0}}\prod_{j=1}^{4}(z-z_{j}).
\end{equation}

But direct computer calculations show that for (\ref{v(z)}) it
holds
$$
v'''-4rv-2r'v=\frac{p\widetilde{P}(z)}{(z-z_{0})^{3}(z^{2}-\varkappa^{2})(z^{2}-\lambda^{2})},
$$
where $\widetilde{P}(z)$ is a polynomial with the leading term
$12i(\varkappa^{2}-\lambda^{2})z^{6}=-48\mu z^{6}/\alpha\ne 0$ for
the spherical case and $12(\lambda^{2}-\varkappa^{2})z^{6}=48\mu
z^{6}/\alpha\ne 0$ for the hyperbolic case. Hence, equation
(\ref{Diff2Equation3}) can not hold.

Thus, in the case I of lemma \ref{Lem1} there can be only one
linear independent solution $y_{1}(z)$ of (\ref{Diff2Equation2})
such that $y_{1}'/y_{1}\in\mathbb{C}(z)$. Since exponents at the
points $z_{j},\,j=1,2,3,4$ are not real one can conclude from
lemma \ref{Lem2} that in this case the group $\EuScript{G}$ is
conjugate to the full triangular group, coincides with
$\EuScript{G}_{0}$ and is not Abelian.

Now, using the Kovacic algorithm, we shall show that the second
case of lemma \ref{Lem1} can not occur. Clearly, it holds $\ord
z_{j}=\ord\infty=2,\,j=0,\ldots,4$ and
$E_{0}=E_{\infty}=(-2,2,6),\,E_{j}=(2),\,j=1,\ldots,4$.\footnote{For
brevity we use the notation $E_{j}:=E_{z_{j}}$.} Therefore one
gets $d(e)=\dfrac12(e_{\infty}-e_{0})-4$ and the maximal value for
$d(e)$ is $0$, which corresponds only to $e_{0}=-2,e_{\infty}=6$.

Thus one should define
$$
\Theta(z):=-\frac1{z-z_{0}}+\sum_{j=1}^{4}\frac1{z-z_{j}}
$$
and verify the equality
\begin{equation}\label{Xi}
\Xi(z):=\Theta''+3\Theta\Theta'+\Theta^{3} -4r\Theta-2r'=0.
\end{equation}
But computer calculations shows that
$$
\Xi(z)=\frac{pP_{*}(z)}{\prod_{j=0}^{4}(z-z_{j})^{2}},
$$
where  $P_{*}(z)$ is a polynomial with the leading term
$3i(\varkappa^{2}-\lambda^{2})z^{6}$ for the spherical case and
$3(\lambda^{2}-\varkappa^{2})z^{6}$ for the hyperbolic case. Thus
the second case of lemma \ref{Lem1} can not occur.

One can conclude from lemma \ref{Lem1} that the differential
Galois group for equation (\ref{Diff2Equation2}) with $r(z)$ given
by (\ref{r(z)}), (\ref{r(z)SN}) or (\ref{r(z)}), (\ref{r(z)HN}) is
either full triangular group (\ref{Tgroup}) or
$\SL_{2}(\mathbb{C})$. In both cases it coincides with its
identity component $\EuScript{G}_{0}$ and is not Abelian.
\end{proof}

Note that due to (\ref{rForMu=1}) for $\mu=1$ equation
(\ref{Diff2Equation2}) with $r(z)$ given by (\ref{r(z)}),
(\ref{r(z)SN}) or (\ref{r(z)}), (\ref{r(z)HN}) has linear
independent solutions $y_{1}=(z-z_{0})^{3/2}$ and
$y_{2}=(z-z_{0})^{-1/2}$. Therefore, for $\mu=1$ the second case
of lemma \ref{Lem2} occurs and $\EuScript{G}=\mathbb{Z}_{2}$.
Here, $\EuScript{G}_{0}$ is trivial that corresponds to the
existence of the additional integral
$p_{1}\sin\theta+p_{2}\cos\theta$ for Hamiltonian function
(\ref{hsSimple}) and $p_{1}\sinh\theta+p_{2}\cosh\theta$ for
Hamiltonian function (\ref{hhSimple}).

\begin{Lem}\label{L6}
Suppose that assumptions of lemma \ref{L3} are valid. Then the
identity component $\EuScript{G}_{0}$ of the Galois group for
equation (\ref{Diff2Equation2}) with $r(z)$ given by
(\ref{r(z)Osc}), (\ref{coeffr(z)SO}) or (\ref{r(z)Osc}),
(\ref{r(z)HO}) is not Abelian.
\end{Lem}
\begin{proof}
There are five regular singular points $z_{j},\,j=0,\ldots,4$, of
order $2$ and the regular singular point $z_{5}=\infty$ of order
$0$. One has
$\Delta_{0}=2,\,\Delta_{j}=\sqrt{1+4\alpha_{j}},\,j=1,2,\,\Delta_{3,4}=1/2$
and due to lemma \ref{L3} $\Delta_{j}\not\in\mathbb{R},\,j=1,2$.
Therefore the third case from lemma \ref{Lem1} is impossible.

Consider the first case of lemma \ref{Lem1}. Suppose that there
are two linear independent solutions $y_{k}(z),\,k=1,2$ of
(\ref{Diff2Equation2}) such that $y_{k}'/y_{k}\in\mathbb{C}(z)$.
Reasoning as in the proof of lemma \ref{L5} one can write them in
the form (\ref{F}), where
$$\rho_{k}^{(0)}\in\left(-\frac12,\frac32\right),\,
\rho_{k}^{(j)}\in\left(\frac12(1+\Delta_{j}),\frac12(1-\Delta_{j})\right),
\,j=1,2,\;\rho_{k}^{(j)}\in\left(\frac14,\frac34\right),\,j=3,4.$$
Besides it holds $v(z):=y_{1}(z)y_{2}(z)\in\mathbb{C}(z)$.

Hence, possible exponents for $v(z)$ are $-1,1,3$ at the point
$z_{0}$; $1,1\pm\sqrt{1+4\alpha_{j}}\not\in\mathbb{R}$ at points
$z_{j},\,j=1,2$; and $1/2,1,3/2$ at points $z_{j},\,j=3,4$. The
inclusion $v(z)\in\mathbb{C}(z)$ implies therefore only two
possibilities:
$$
v(z)=P(z)\prod_{j=1}^{4}(z-z_{j})\quad\text{or}\quad
v(z)=\frac{P(z)}{z-z_{0}}\prod_{j=1}^{4}(z-z_{j}),\;P(z)\in\mathbb{C}[z].
$$
Due to (\ref{r(z)Osc}) exponents for (\ref{Diff2Equation2}) at
$\infty$ are $0$ and $-1$, the function $v(z)$ grows as
$z\to\infty$ no faster than $z^{2}$; therefore no one of these
possibilities can realize.

Thus in the case I of lemma \ref{Lem1} there can be only one
linear independent solution $y_{1}(z)$ of (\ref{Diff2Equation2})
such that $y_{1}'/y_{1}\in\mathbb{C}(z)$. Reasoning as in the
proof of lemma \ref{L5} one gets that in this case the group
$\EuScript{G}_{0}$ is conjugate to the full triangular group and
is not Abelian.

Check the possibility of the second case of lemma \ref{Lem1} using
the Kovacic algorithm. Clearly, it holds $\ord
z_{j}=2,\,j=0,\ldots,4,\;\ord\infty=0$ and
$E_{0}=(-2,2,6),\,E_{1,2}=(2),\,E_{3,4}=(1,2,3),\,E_{\infty}=(0,2,4)$.
Therefore, a unique element $e\in E$ for which $d(e)\geqslant0$ is
$(-2,2,2,1,1,4)$ and $d(-2,2,2,1,1,4)=0$.

Thus one should define
$$
\Theta(z):=-\frac1{z-z_{0}}+\frac1{z-z_{1}}+\frac1{z-z_{2}}+\frac1{2(z-z_{3})}
+\frac1{2(z-z_{4})}
$$
and verify equality (\ref{Xi}).

But computer calculations shows that
$$
\Xi(z)=\frac{2p(2\mu+1)(\varkappa^{2}-\lambda^{2})^{2}z^{2}-8p^{2}(\mu-1)
(\varkappa^{2}-\lambda^{2})^{3}z+2p^{2}(\mu-1)(\varkappa^{2}-\lambda^{2})^{2}-
3\varkappa^{2}}{(z-z_{3})(z-z_{4})\prod_{j=0}^{2}(z-z_{j})^{2}}.
$$
for the spherical case and
$$
\Xi(z)=\frac{-2p(2\mu+1)(\varkappa^{2}-\lambda^{2})^{2}z^{2}-8p^{2}(\mu-1)
(\varkappa^{2}-\lambda^{2})^{3}z-2p^{2}(\mu-1)(\varkappa^{2}-\lambda^{2})^{2}+
3\varkappa^{2}}{(z-z_{3})(z-z_{4})\prod_{j=0}^{2}(z-z_{j})^{2}}.
$$
for the hyperbolic case. Thus the second case of lemma \ref{Lem1}
can not occur.

Now lemma \ref{Lem1} implies that the differential Galois group
for equation (\ref{Diff2Equation2}) with $r(z)$ given by
(\ref{r(z)Osc}), (\ref{coeffr(z)SO}) or (\ref{r(z)Osc}),
(\ref{r(z)HO}) is either full triangular group (\ref{Tgroup}) or
$\SL_{2}(\mathbb{C})$. In both cases it coincides with its
identity component $\EuScript{G}_{0}$ and is not Abelian.
\end{proof}

Note that due to (\ref{rOscMu=1}) for $\mu=1$ equation
(\ref{Diff2Equation2}) with $r(z)$ given by (\ref{r(z)Osc}),
(\ref{coeffr(z)SO}) or (\ref{r(z)Osc}), (\ref{r(z)HO}) has linear
independent solutions
$$
y_{1}(z)=\frac{\left(z^{2}-\lambda^{2}\right)^{3/4}}{\left(z-z_{0}\right)^{1/2}}
\quad\text{and}\quad
y_{2}(z)=\frac{\left(z^{2}-\lambda^{2}\right)^{1/4}\left(z_{0}z-\lambda^{2}\right)}
{\left(z-z_{0}\right)^{1/2}}.
$$ Therefore for $\mu=1$ the second case of
lemma \ref{Lem2} occurs, $\EuScript{G}=\mathbb{Z}_{4}$ and
$\EuScript{G}_{0}$ is trivial. This corresponds to the existence
of the additional integrals described above.

{\it Proof of theorems \ref{MainThS} and \ref{MainThH}}. Due to
theorem \ref{MRTh} and the analysis in section \ref{ParSol} it is
enough to show that the identity components of Galois groups for
systems (\ref{NVESz}) and (\ref{NVESzOs}) described in section
\ref{ParSol} are not Abelian.

Consider transformations of system (\ref{NVESz}) made in section
\ref{ParSol1}. The base field for (\ref{NVESz}) is
$\mathbb{C}(z)$. First, we reduced (\ref{NVESz}) to linear
differential equation (\ref{p21}) of the second order that
corresponds to the variable change
$$
\left(\begin{array}{c} p_{1} \\ p_{2} \end{array}\right)\to
\left(\begin{array}{c} p_{2}' \\ p_{2} \end{array}\right)=
\left(\begin{array}{cc} C(z) & -A(z) \\ 0 & 1 \end{array}\right)
\left(\begin{array}{c} p_{1} \\ p_{2} \end{array}\right),
$$
which is reversible over $\mathbb{C}(z)$. Therefore Galois groups
for system (\ref{NVESz}) and equation (\ref{p21}) coincides. Then
we came to equation (\ref{Diff2Equation2}) for the function
$y(z)=p_{2}(z)/\sqrt{C(z)}$ with $r(z)$ given by (\ref{r(z)}).
Since the function $\sqrt{C(z)}$ is algebraic the identity
components of Galois groups for equations (\ref{p21}) and
(\ref{Diff2Equation2}) are the same. One completes the proof for
system (\ref{NVESz}) using lemmas \ref{L1} and \ref{L2}.

The proof of theorem \ref{MainThH} is similar due to lemma
\ref{L3} and since all transformations made in section
\ref{ParSol2} from system (\ref{NVESzOs}) to equation
(\ref{Diff2Equation2}) are linear and their coefficients are
algebraic functions.\hfill$\square$

\section{Conclusion}\label{Conclusion}\markright{\ref{Conclusion}
Conclusion}

In the present paper we proved the complex nonintegrability of the
reduced two-body problem in the spaces ${\bf S}^{2}$ and ${\bf
H}^{2}$ for the Newton and the oscillator potentials. The main
prerequisite for this proof was the possibility to reduce the
system of normal variational equations to a linear differential
equation of the second order with rational coefficients using a
proper change of an independent variable. It is obvious that for a
general central potential it could not be done.

Therefore, the problem of finding a nontrivial central potential
corresponding to integrability of the two-body problem in constant
curvature spaces or proving the absence of such potential in some
more or less general class is open.

\appendix
\section{Appendix}\label{Appendix}\markboth{\ref{Appendix} Appendix}
{\ref{Appendix} Appendix}

We use the standard notations $\mathbb{C}(z)$ and $\mathbb{C}[z]$
for the field of rational functions and for the ring of
polynomials, both with complex coefficients. Consider a linear
second order differential equation on the Riemannian sphere
$\mathbf{P}^{1}(\mathbb{C})$
\begin{equation}\label{Diff2Equation}
w''(z)+p(z)w'(z)+q(z)w(z)=0,\;p(z),q(z)\in\mathbb{C}(z).
\end{equation}
Any pole $z_{0}\in\mathbb{C}$ of $p(z)$ or $q(z)$ is a singular
point of equation (\ref{Diff2Equation}). This point is a {\it
regular singular point} for (\ref{Diff2Equation}) if functions
$(z-z_{0})p(z)$ and $(z-z_{0})^{2}q(z)$ are holomorphic at
$z_{0}$.  One can find {\it exponents} $\rho^{(z_{0})}$ of
(\ref{Diff2Equation}) at the point $z_{0}$ by the substitution
$w(z)=(z-z_{0})^{\rho^{(z_{0})}}$ into (\ref{Diff2Equation}) and
keeping only leading terms as $z\rightarrow z_{0}$. This procedure
gives a quadratic equation for $\rho^{(z_{0})}$.

The same is also valid for the point
$\infty\in\mathbf{P}^{1}(\mathbb{C})$ w.r.t.\ a variable
$\zeta=1/z$.

Equation (\ref{Diff2Equation}) is {\it Fuchsian} iff all its
singular points are regular.

\begin{Def}
A solution $\widetilde{w}(z)$ of equation (\ref{Diff2Equation}) is
called Liouvillian if there is a tower of differential fields
$$
\mathbb{C}(z)=\mathcal{K}_{0}\subset\mathcal{K}_{1}\subset\ldots\subset\mathcal{K}_{m}
$$
with $\widetilde{w}(z)\in\mathcal{K}_{m}$ and for each
$i=1,\ldots,m$ it holds
$\mathcal{K}_{i}=\mathcal{K}_{i-1}(v_{i})$, where one of the
following three possibilities holds:
\begin{enumerate}
\item the function $v_{i}$ is algebraic over $\mathcal{K}_{i-1}$;
\item it satisfies to $v_{i}'\in\mathcal{K}_{i-1}$;
\item it satisfies to $v_{i}'/v_{i}\in\mathcal{K}_{i-1}$.
\end{enumerate}
\end{Def}

The substitution
$$
w(z)=\exp\left(-\frac12\int p(z)dz\right)y(z)
$$
transforms equation (\ref{Diff2Equation}) into the equation
\begin{equation}\label{Diff2Equation2}
y''(z)=r(z)y(z),
\end{equation}
where $r(z):=-q(z)+\frac12p'(z)+\frac14p^{2}(z)$. If
$y_{1}(z),y_{2}(z)$ are solutions of (\ref{Diff2Equation}), then a
direct calculation shows that
\begin{equation}\label{Diff2Equation3}
v'''(z)=4rv'(z)+2r'v(z)
\end{equation}
for $v(z):=y_{1}(z)y_{2}(z)$. Equation (\ref{Diff2Equation3}) is
called the second symmetric power of (\ref{Diff2Equation2}).

The differential Galois group $\EuScript{G}$ for
(\ref{Diff2Equation2}) is an algebraic subgroup of
$\SL_{2}(\mathbb{C})$ \cite{Kap}. The following lemma \cite{Kov}
contains a classification of possible groups $\EuScript{G}$.
\begin{Lem}\label{Lem1}
One and only one of the following four cases can occur.\\
\vskip1pt\noindent Case I. The group $\EuScript{G}$ is conjugate
to a subgroup of the full triangular group
\begin{equation}\label{Tgroup}
T=\left\{\left.\left(\begin{array}{cc} a & b \\ 0 & a^{-1}
\end{array}\right)\right|\;a\in\mathbb{C}^{*},b\in\mathbb{C}\right\};
\end{equation}
in this case equation (\ref{Diff2Equation2}) has a solution
$y_{1}(z)\not\equiv0$ such that $y'_{1}/y_{1}\in\mathbb{C}(z)$.\\
\vskip1pt\noindent Case II. The group $\EuScript{G}$ is conjugate
to a subgroup of
$$
\left\{\left.\left(\begin{array}{cc} c & 0 \\ 0 & c^{-1}
\end{array}\right)\right|\;c\in\mathbb{C}^{*}\right\}\cup
\left\{\left.\left(\begin{array}{cc} 0 & c \\ -c^{-1} & 0
\end{array}\right)\right|\;c\in\mathbb{C}^{*}\right\}
$$
and case I does not hold. In this case  equation
(\ref{Diff2Equation2}) has a solution of the form
$y_{1}(z)=\exp\left(\int\omega(z)dz\right)$, where $\omega(z)$ is
an algebraic function over $\mathbb{C}(z)$ of degree 2.\\
\vskip1pt\noindent Case III. The group $\EuScript{G}$ is finite
and cases I and II do not hold. In this case all solutions of
(\ref{Diff2Equation2}) are algebraic and exponents of
(\ref{Diff2Equation2}) are rational numbers at all points. \\
\vskip1pt\noindent Case IV.
$\EuScript{G}\simeq\SL_{2}(\mathbb{C})$ and equation
(\ref{Diff2Equation2}) has no Liouvillian solutions.
\end{Lem}

A more precise information on the case I from the preceding lemma
is contained in the following lemma.
\begin{Lem}[see proposition 4.2 from \cite{SU}]\label{Lem2}
Suppose that the case I from lemma \ref{Lem1} occurs and therefore
$\EuScript{G}$ is conjugate to a subgroup of the group $T$.
\begin{enumerate}
\item If equation (\ref{Diff2Equation2}) has a unique (up to a constant
factor) solution $y_{1}(z)\not\equiv0$ such that
$y'_{1}/y_{1}\in\mathbb{C}(z)$, then $\EuScript{G}$ is conjugate
to a proper subgroup of the group $T$ iff
$y_{1}^{m}\in\mathbb{C}(z)$ for some $m\in\mathbb{N}$. In this
case $\EuScript{G}$ is conjugate to
$$
T_{m}=\left\{\left.\left(\begin{array}{cc} a & b \\ 0 & a^{-1}
\end{array}\right)\right|\;a,b\in\mathbb{C},a^{m}=1\right\},
$$
where $m$ is the smallest positive integer such that
$y_{1}^{m}\in\mathbb{C}(z)$.
\item If equation (\ref{Diff2Equation2}) has two linear independent
solutions $y_{1},y_{2}$ such that
$y'_{j}/y_{j}\in\mathbb{C}(z),\,j=1,2$, then $\EuScript{G}$ is
conjugate to a subgroup of the group
$$
D=\left\{\left.\left(\begin{array}{cc} a & 0 \\ 0 & a^{-1}
\end{array}\right)\right|\;a\in\mathbb{C}^{*}\right\}.
$$
In this case $y_{1}y_{2}\in\mathbb{C}(z)$. Finally, $\EuScript{G}$
is conjugate to a proper subgroup of the group $D$ iff
$y_{1}^{m}\in\mathbb{C}(z)$ for some $m\in\mathbb{N}$. In this
case $\EuScript{G}$ is conjugate to a cyclic group of order $m$,
where $m$ is the smallest positive integer such that
$y_{1}^{m}\in\mathbb{C}(z)$.
\end{enumerate}
\end{Lem}

We also need the Kovacic algorithm for case II from lemma
\ref{Lem1}. It allows one to find a solution of equation
(\ref{Diff2Equation2}) of the form
$\exp\left(\int\omega(z)dz\right)$, where $\omega(z)$ is an
algebraic function of degree $2$, or conclude that such solution
does not exist.

Let $r(z)=s(z)/t(z)$, where $s(z)$ and $t(z)$ are relatively prime
polynomials and $t(z)$ is monic, i.e.\ the coefficient of its
leading term equals $1$. We denote
$\Sigma':=\left(\left.c\in\mathbb{C}\right|t(c)=0\right)$ and
$\Sigma:=\Sigma'\cup(\infty)$.

Let an order $\ord c$ of $c\in\Sigma'$ be the multiplicity of $c$
as a root of $t(z)$ and $\ord \infty:=\max\left(0,4+\deg s-\deg
t\right)$. If $c\in\Sigma$ such that $\ord c=1$ or $\ord c=2$, one
can find an expansion
$$
r(z)=\frac{\alpha_{c}}{(z-c)^{2}}+O\left(\frac1{z-c}\right)\quad\text{as}\quad
z\to c\quad\text{for}\quad c\in\Sigma'
$$
and
$$
r(z)=\frac{\alpha_{\infty}}{z^{2}}+O\left(\frac1{z^{3}}\right)\quad\text{as}\quad
z\to \infty\quad\text{for}\quad c=\infty.
$$
It these cases $\Delta_{c}:=\sqrt{1+4\alpha_{c}}$ is the
difference of exponents for (\ref{Diff2Equation2}) at $z=c$.
Evidently $\alpha_{c}=0$ and $\Delta_{c}=1$ if $\ord c=1$.

{\bf Step 1}. For every $c\in\Sigma$ we define a finite set
$E_{c}$ in the following way.

If $\ord \infty=0$, then put $E_{\infty}=(0,2,4)$.

If $\ord c=1$, then put $E_{c}=(4)$ for $c\neq\infty$ and
$E_{\infty}=(0,2,4)$.

If $\ord c=2$, then put
$$E_{c}=(2,2(1+\Delta_{c}),2(1-\Delta_{c}))\cap\mathbb{Z}.$$

If $\ord c=k>2$, then put $E_{c}=(k)$ for $c\neq\infty$ and
$E_{\infty}=(4-k)$.

{\bf Step 2}. For each element $e$ of the set
$$
E:=\prod_{c\in\Sigma}E_{c}
$$
we compute the number
$$d(e):=\dfrac12\left(e_{\infty}-\sum_{c\in\Sigma'}e_{c}\right),$$
where $e_{c}$ is a component of an array $e$ from $E_{c}$. We
select those elements $e\in E$ for which
$d(e)\in\mathbb{Z}_{+}:=\mathbb{N}\cup (0)$. If there are no such
elements, then the case II from lemma \ref{Lem1} can not occur.

{\bf Step 3}. For each element $e\in E$ selected on the previous
step we define
$$
\Theta(z):=\frac12\sum_{c\in\Sigma'}\frac{e_{c}}{z-c}
$$
and search for a monic polynomial $P(z)$ of degree $d(e)$
satisfying the following equation
$$
P'''+3\Theta
P''+\left(3\Theta^{2}+3\Theta'-4r\right)P'+\left(\Theta''+3\Theta\Theta'+\Theta^{3}
-4r\Theta-2r'\right)P=0.
$$
If such a polynomial exists, then equation (\ref{Diff2Equation2})
has a solution of the form $\exp\left(\int\omega(z)dz\right)$,
where
$$
\omega^{2}-\psi\omega+\frac12\psi'+\frac12\psi^{2}-r=0,\;\psi=\psi(z):=\Theta+\frac{P'}P.
$$
If such polynomial $P$ does not exist, then case II from lemma
\ref{Lem1} does not occur.

\end{document}